\documentclass[12pt,dvips]{article}
\usepackage{amsmath,amsfonts,amssymb,amsthm,graphicx,verbatim,subfig}
\usepackage[all]{xy}

\ifx\pdftexversion\undefined

\usepackage[a4paper,colorlinks,link
=black,filecolor=black,citecolor=black,urlcolor=black,pdfstartview=FitH]{hyperref}
\else

\usepackage[a4paper,colorlinks,linkcolor=black,filecolor=black,citecolor=black,urlcolor=black,pdfstartview=FitH]{hyperref}
\fi

\font\sixbb=msbm6
\font\eightbb=msbm8
\font\twelvebb=msbm10 scaled 1095
\newfam\bbfam
\textfont\bbfam=\twelvebb \scriptfont\bbfam=\eightbb
                           \scriptscriptfont\bbfam=\sixbb
\def\bb{\fam\bbfam\twelvebb}
\newcommand{\Rea}{{\bb R}}

\newcommand{\FF}{{\bb F}}

\newtheorem{theorem}{\bf Theorem}
\newtheorem{claim}[theorem]{\bf Claim}

\newtheorem{proposition}[theorem]{\bf Proposition}
\newtheorem{corollary}[theorem]{\bf Corollary}
\newtheorem*{definition}{Definition}
\newcommand{\enp}{\begin{flushright} $\Box$ \end{flushright}}
\newcommand{\beq}[0]{\begin{equation}}
\newcommand{\enq}[0]{\end{equation}}

\newcommand{\dn}{\Delta_{n-1}}

\newcommand{\cg}{{\cal G}}

\newcommand{\supp}{{\rm supp}}

\newcommand{\cl}{{\cal L}}

\newcommand{\pr}{{\rm Pr}}

\newcommand{\lk}{\text{lk}}
\newcommand{\per}{\text{per}}
\newcommand{\cy}{{\cal Y}}
\newcommand{\ce}{{\cal E}}
\newcommand{\ccs}{{\mathbb S}}
\newcommand{\ul}{\underline{L}}

\title{Random Latin Squares and \\ $2$-Dimensional Expanders}
\begin{document}
\author{Alexander Lubotzky\thanks{Institute of Mathematics, Hebrew University, Jerusalem 91904, Israel. e-mail: alexlub@math.huji.ac.il~. Supported by ERC and NSF grants.} and Roy Meshulam\thanks{Department of Mathematics,
Technion, Haifa 32000, Israel. e-mail:
meshulam@math.technion.ac.il~. Supported by an ISF grant.}}
\maketitle
\pagestyle{plain}
\begin{abstract}
Let $X$ be a $2$-dimensional simplicial complex. The degree of an edge $e$ is the number of $2$-faces of $X$ containing $e$.
The complex $X$ is an $\epsilon$-expander if the coboundary $d_1\phi$ of every $\FF_2$-valued $1$-cochain $\phi \in C^1(X;\FF_2)$ satisfies $$|\supp(d_1\phi)| \geq \epsilon \min\{|\supp(\phi+d_0\psi)|: \psi \in C^0(X;\FF_2)\}.$$

In response to \cite{G10} and \cite{DK12} we show the existence of $2$-dimensional $\epsilon$-expanders with maximum edge degree $d$ for some fixed $\epsilon>0$ and $d$.
This is done via the following new model of random $2$-dimensional complexes.
A Latin square of order $n$ is an $n$-tuple $L=(\pi_1,\ldots,\pi_n)$ of permutations on $\{1,\ldots,n\}$ such that
$\pi_i^{-1} \pi_j$ is fixed point free for $1 \leq i<j \leq n$. Let $\{a_i\}_{i=1}^n, \{b_i\}_{i=1}^n, \{c_i\}_{i=1}^n$ be three disjoint sets and
let $(L_1,\ldots,L_d)$ be a $d$-tuple of independently chosen random Latin squares
of order $n$. For $1 \leq k \leq d$, let  $L_k=(\pi_{k1},\ldots,\pi_{kn})$. Let $Y(L_1,\ldots,L_d)$ be the $3$-partite $2$-dimensional complex consisting of all
$2$-simplices $[a_i,b_j,c_{\pi_{ki}(j)}]$ where $1 \leq i,j \leq n$ and $1 \leq k \leq d$.

It is shown that there exist $d< \infty$ and $\epsilon>0$ such that the complex $Y(L_1,\ldots,L_d)$ is an $\epsilon$-expander
with probability tending to 1 as $n \rightarrow \infty$.
\end{abstract}

\section{Introduction}
\label{intro}
\ \ \
The notion of expander graphs has been extremely useful in computer science, combinatorics and even pure mathematics (see \cite{HLW06,Lu2} and the references therein).
In recent years there is a growing interest in high-dimensional expanders (see the survey \cite{Lu3}).
While there are several "competing" definitions for $k$-dimensional expansion, we shall focus
on the notion of {\it coboundary expansion} of simplicial complexes. This  $k$-dimensional version of the graphical Cheeger constant came up independently in the work of Linial, Meshulam and Wallach \cite{LM06,MW09} on homological connectivity of random complexes and in Gromov's remarkable work \cite{G10,FGLNP} where it is shown that this notion of expansion implies the topological overlap property (see Section \ref{crem}).

We recall some topological terminology.
Let $X$ be a simplicial complex on the vertex set $V$. For $k \geq 0$, let $X^{(k)}$ denote the $k$-dimensional skeleton of $X$ and let
$X(k)$ be the family of $k$-dimensional faces of $X$. Let $D_k(X)$ be the maximum number of $(k+1)$-dimensional faces of $X$ containing a common $k$-dimensional face.
Let $C^k(X;\FF_2)$ denote the space of  $\FF_2$-valued functions on $X(k)$.
The $k$-coboundary map $d_k:C^k(X;\FF_2) \rightarrow C^{k+1}(X;\FF_2)$ is given
$$d_k\phi(v_0,\cdots,v_{k+1})=\sum_{i=0}^{k+1} \phi(v_0,\ldots,v_{i-1}, v_{i+1}, \ldots,v_{k+1}).$$
It will be convenient to augment the cochain
complex $\{C^i(X;\FF_2)\}_{i=0}^{\infty}$ with a $(-1)$-degree term
$C^{-1}(X;\FF_2)=\FF_2$ with a coboundary map $d_{-1}:C^{-1}(X;\FF_2)
\rightarrow C^0(X;\FF_2)$ given by $d_{-1}(a)(v)=a$ for $a \in \FF_2, v
\in V$.

Let $Z^k(X;\FF_2)=\ker d_k$ be the space of  $k$-dimensional $\FF_2$-cocycles of $X$ and
let $B^k(X;\FF_2)=d_{k-1}(C^{k-1}(X;\FF_2))$ be the space of  $k$-dimensional $\FF_2$-coboundaries of $X$.
The $k$-dimensional reduced $\FF_2$-cohomology group of $X$ is $$\tilde{H}^k(X;\FF_2)=\frac{Z^k(X;\FF_2)}{B^k(X;\FF_2)}.$$
For $\phi \in C^k(X;\FF_2)$, let $[\phi]$ denote the image of $\phi$ in the quotient space
$C^k(X;\FF_2)/B^k(X;\FF_2)$.
Let $$\|\phi\|=|\{\sigma \in X(k): \phi(\sigma) \ne 0\}|=|\supp(\phi)|$$
and  $$\|[\phi]\|=\min\{|\supp(\phi+d_{k-1}\psi)|: \psi \in C^{k-1}(X;\FF_2)\}.$$
We will sometimes write $\|\phi\|_X$ and $\|[\phi]\|_X$ in case of ambiguity concerning $X$.
\begin{definition}
The {\it $k$-th coboundary expansion constant} of $X$ (see \cite{LM06,MW09,G10,DK12}) is defined by
$$
h_k(X)=\min \left\{\frac{\|d_k\phi \|}{\|[\phi]\|}: \phi \in C^k(X;\FF_2)-B^k(X;\FF_2)\right\}.
$$
A complex $X$ is a {\it $(k,d,\epsilon)$-expander} if $$D_{k-1}(X) \leq d  ~~~~\text{and}~~~~  h_{k-1}(X) \geq \epsilon.$$
\end{definition}
\noindent
{\bf Remarks:}\\
1. Note that $h_k(X)=0$ iff $\tilde{H}^k(X;\FF_2) \neq 0$.
\\
2. Let $G=(V,E)$ be a graph. For $S \subset V$ let $$E(S,V-S)=\{e \in E:|e \cap S|=1\}$$ be the cut determined by $S$ and let $e(S,V-S)=|E(S,V-S)|$.
Viewing $G$ as a $1$-dimensional complex, it is easy to check (see \cite{DK12}) that $h_0(G)$ coincides with the Cheeger constant of $G$:
$$h_0(G)=\min_{0< |S| \leq \frac{|V|}{2}} \frac{e(S,V-S)}{|S|}.$$
\\
3. Expansion in graphs can be defined in two essentially equivalent ways, either via the Cheeger constant or via the spectral gap of the graph Laplacian. This equivalence does not however hold in higher dimensions. Indeed, while spectral gaps of the $k$-dimensional Laplacians carry substantial information concerning the combinatorics and topology of a complex, they cannot in general detect homology in positive characteristics.
\vspace{3mm}

Questions concerning existence and construction of families of expanders graphs have been the subject of intense research for the last 40 years.
Clearly the complete graph is an expander. Using random methods it is also not difficult to show the existence of infinite
families of $(1,d,\epsilon)$-expander graphs for some fixed $d$ and $\epsilon>0$.
Explicit constructions are considerably deeper and involve applications of Kazhdan Property (T), the Ramanujan conjecture and the Zigzag product (see \cite{Lu1,HLW06}).

In this paper we are concerned with the existence of higher dimensional expanders.
Here even the fact that the full $k$-skeleton of the
$(n-1)$-simplex $\dn$ is expanding is not completely obvious.
In \cite{MW09,G10} it was shown that the expansion of $\dn$ is given by
\begin{equation}
\label{expsim}
h_{k-1}(\dn)=\frac{n}{k+1}.
\end{equation}

We next consider the expansion of random complexes in the following natural $k$-dimensional extension of the
Erd\H{o}s-R\'{e}nyi random graph model (see \cite{LM06,MW09}).
For $k \geq 1$ and $0 \leq p \leq 1$ let
$Y_k(n,p)$ denote the probability space of all $k$-dimensional random subcomplexes of $\dn$ obtained
by starting with the full $(k-1)$-dimensional skeleton of $\dn$ and then adding
each $k$-simplex independently with probability $p$.
Using the Chernoff bound it directly follows from (\ref{expsim})
that there exists a constant $C(k)=\Theta(k^2)$ such that if $p=\frac{C(k) \log n}{n}$ then $Y\in Y_k(n,p)$ is asymptotically almost surely (a.a.s.) a $(k,\Theta(\log n),1)$-expander (see \cite{GW12,DK12}).  This of course implies that $Y\in Y_k(n,p)$ a.a.s. satisfies $H_{k-1}(Y;\FF_2)=0$. Obtaining the exact probability threshold $p=\frac{k \log n}{n}$ for the vanishing of $H_{k-1}(Y;\FF_2)$ is more involved (see \cite{LM06,MW09}).
Dotterrer and Kahle \cite{DK12} proved results analogous to (\ref{expsim}) for skeleta
of cross polytopes and for complete multipartite complexes with similar consequences for the expansion of their random subcomplexes. As in the case of $\dn$, all these complexes have unbounded degrees.

Here we establish the existence of an infinite family of
$(2,d,\epsilon)$-expanders for some fixed $d$ and $\epsilon>0$,
answering a question asked implicitly in \cite{G10} and explicitly in \cite{DK12}. Our proof is probabilistic and depends
on the following new model, based on Latin squares, of random $2$-dimensional simplicial complexes with bounded edge degrees.

We introduce some notation and definitions that will be used throughout the paper.
Let $\ccs_n$ be the symmetric group on $[n]=\{1,\ldots,n\}$.
A $k$-tuple $(\pi_1,\ldots,\pi_k) \in \ccs_n^k$ is {\it legal} if $\pi_i\pi_j^{-1}$ is fixed point free for all $1 \leq i < j \leq k$.
A {\it Latin Square} of order $n$ is a legal $n$-tuple of permutations $L=(\pi_1,\ldots,\pi_n) \in \ccs_n^n$.
Let $\cl_n$ denote the uniform probability space of all Latin squares of order $n$.
Let $V_1=\{a_i\}_{i=1}^n, V_2=\{b_i\}_{i=1}^n, V_3=\{c_i\}_{i=1}^n$ be three disjoint sets.
The complete $3$-partite complex $T_n=V_1*V_2*V_3$ consists of all $\sigma \subset V=V_1 \cup V_2 \cup V_3$
such that $|\sigma \cap V_i| \leq 1$ for $1 \leq i \leq 3$.
An $L=(\pi_1,\ldots,\pi_n) \in \cl_n$ determines a subcomplex $T_n^{(1)} \subset Y(L) \subset T_n$ whose $2$-simplices are
$[a_i,b_j,c_{\pi_i(j)}]$ where $1 \leq i,j \leq n$.
Fix $d$ and regard  $\cl_n^d$ as a uniform probability space.
For $\ul^d=(L_1,\ldots,L_d) \in \cl_n^d$, let $Y(\ul^d)= \cup_{i=1}^d Y(L_i)$.
Note that $D_1(Y(\ul^d)) \leq d$.
Let $\cy(n,d)$ denote the probability space of all complexes $Y(\ul^d)$ with measure induced from $\cl_n^d$.
\begin{theorem}
\label{main}
There exist $\epsilon>0, d<\infty$ such that
$$\lim_{n \rightarrow \infty} \pr\left[Y \in \cy(n,d): h_1(Y)>\epsilon\right]=1.$$
\end{theorem}
Theorem  \ref{main} is proved in two steps: in the first it is shown that $1$-cochains with small support have large coboundaries. This part is inspired by (and is in fact easier than) the results of \cite{KKL13}. To handle the case of $1$-cochains with large support we have to establish a concentration of measure theorem for the space $\cl_n$.
This result may be viewed as a (coarse) Latin square analogue of Maurey's large deviations bound for Lipschitz functions on the symmetric group \cite{Ma79} (for further comments see Section \ref{s:ldls}).

The paper is organized as follows. In Section \ref{s:outline} we describe the strategy of the proof. In Section \ref{spect} we prove a lower bound (Proposition \ref{onetwo1}) on the expansion of small cochains in terms of the spectral gaps of the vertex links. In Section \ref{s:ldls} we prove a large deviations bound for random Latin squares (Theorem \ref{largedev}). This result is the key for showing expansion of large cochains. In Section \ref{prf} we use the above mentioned results to obtain Theorem \ref{main}. We conclude in Section \ref{crem} with some questions and comments.

\section{Outline of the proof}
\label{s:outline}
\ \ \
In order to prove Theorem \ref{main} we have to show the existence of fixed $d$ and $\epsilon>0$ such that for
almost all $Y \in \cy(n,d)$, all $\phi \in C^1(T_n;\FF_2)$ satisfy $\|d_1\phi\|_Y \geq \epsilon \|[\phi]\|$.

Let $0<c<1$ be a constant whose value will be assigned later. A $1$-cochain $\phi$ is called {\it $c$-small} if $\|[\phi]\| \leq cn^2$
and {\it $c$-large} otherwise. We first address the expansion of $c$-small $1$-cochains. Let $Y$ be a subcomplex of $T_n$ with maximum edge degree at most $d$.
The link $Y_v=\lk(Y,v)$ of a vertex $v \in V$ is a bipartite graph. Let $\mu(Y_v)$ be the second smallest eigenvalue of the Laplacian of $Y_v$ and
let $\tilde{\mu}=\min_{v \in V}\mu(Y_v)$. Inspired by the results of \cite{KKL13}, we show in Proposition \ref{onetwo1} that
\begin{equation}
\label{out1}
\|d_1\phi\|_Y \geq \left(\frac{(1-c^{1/3})\tilde{\mu}}{2}-\frac{d}{3}\right) \|[\phi]\|.
\end{equation}
Suppose now that $Y=Y(\ul^d) \in \cy(n,d)$.
As each $Y_v$ is a random bipartite graph of maximum degree $d$, it follows from a theorem of Friedman \cite{F91,F08} that $\mu(Y_v)\geq d-O(\sqrt{d})$
for all $v \in V$ with probability $1-o(1)$. Hence by (\ref{out1}) there exist $c>0$ and $d$ such that for almost all $Y \in \cy(n,d)$
\begin{equation}
\label{out3}
\|d_1\phi\|_Y \geq  \|[\phi]\|
\end{equation}
for all $c$-small $\phi$'s.

We next consider the expansion of $c$-large $1$-cochains. Here, in contrast with the first case, we bound the probability of non-expansion separately for each cochain.
Let $\phi \in C^1(T_n;\FF_2)$ be $c$-large.  Dotterrer and Kahle \cite{DK12} proved that
$h_1(T_n) \geq \frac{n}{5}$. It follows that
$$\ce=\{\sigma \in T_n(2): d_1\phi(\sigma) \neq 0\}$$
satisfies $$|\ce|=\|d_1\phi\|_{T_n} \geq \frac{n}{5} \|[\phi]\| \geq \frac{cn^3}{5}.$$
For a Latin square $L$ in the uniform probability space $\cl_n$ let $f_{\ce}(L)=|Y(L)(2) \cap \ce|$ be the number of $2$-simplices in $Y(L)$ that belong to $\ce$.
The expectation of $f_{\ce}$ over $\cl_n$ is
$$E[f_{\ce}]=\frac{|\ce|}{n} \geq \frac{cn^2}{5}.$$
In Theorem \ref{largedev} we prove a large deviation type bound for $f_{\ce}$ showing that
$$\pr[f_{\ce}<10^{-5}c^2n^2] < e^{-\delta n^2}$$
for some $\delta=\delta(c)>0$.
Now let $\ul^d=(L_1,\ldots,L_d) \in \cl_n^d$.
Since
\begin{equation*}
\begin{split}
&\|d_1 \phi\|_{Y(\ul^d)}=|Y(\ul^d)(2) \cap \ce|  \\
&\geq \max_{1 \leq i \leq d} |Y(L_i)(2) \cap \ce|=\max_{1 \leq i \leq d} f_{\ce}(L_i)
\end{split}
\end{equation*}
it follows that
\begin{equation}
\label{out2}
\pr\left[\ul^d \in \cl_n^d: \|d_1\phi\|_{Y(\ul^d)} < 10^{-5}c^2n^2\right] <e^{-\delta d n^2}.
\end{equation}
As the number of $1$-cochains is $2^{3n^2}$, Eq. (\ref{out2}) implies that the probability that there exists a large $1$-cochain $\phi$  with
$\|d_1\phi\|_{Y(\ul^d)} < 10^{-5}c^2n^2$ is at most $2^{3n^2}e^{-\delta d n^2}$. Choosing a sufficiently large $d$ and noting that $\|[\phi]\| \leq 3n^2$ for all $\phi$, it follows that
\begin{equation}
\label{out4}
\pr\left[\ul^d \in \cl_n^d:\frac{\|d_1\phi\|_{Y(\ul^d)}}{\|[\phi]\|} \geq 10^{-6}c^2 \text{~~for~all~c-large~} \phi~ \right] =1-o(1).
\end{equation}

Theorem \ref{main} is now a consequence of (\ref{out3}) and (\ref{out4}). In the next sections we prove
Proposition \ref{onetwo1} and Theorem \ref{largedev} which are the key ingredients in the argument.

\section{Expansion of Small Cochains}
\label{spect}
\ \ \
In this section we relate the expansion of small $1$-cochains to the spectral gaps of the vertex links.
The Laplacian of a graph $G=(V,E)$ is the $V \times V$ positive semidefinite matrix
$L_G$ given by $$L_G(u,v) = \left\{
\begin{array}{ll}
        \deg(u)  & u=v \\
        -1       & uv \in E \\
         0   & {\rm otherwise}
\end{array}
\right.~~ $$ Let $0=\mu_1(G) \leq \mu_2(G) \leq \cdots
\leq \mu_{|V|}(G)$ denote the eigenvalues of $G$. The spectral gap of $G$ is $\mu(G)=\mu_2(G)$.

Let $Y$ be a subcomplex of the complete $3$-partite complex $T_n$ on the vertex set $V=V_1 \cup V_2 \cup V_3$ with a complete $1$-skeleton $T_n^{(1)}$.
The link  $Y_v=\lk(Y,v)$ of a vertex $v \in V_i$ is a bipartite graph on the vertex set $V_v=\lk(Y,v)(0)=V_j \cup V_{j'}$ where $\{1,2,3\}=\{i,j,j'\}$.
Let $d=\max_{e \in Y(1)} |\lk(Y,e)|$ be the maximal degree of an edge in $Y$.
Let $\tilde{\mu}=\min_{v \in V} \mu(Y_v)$.

Let $\phi \in C^1(Y;\FF_2)$ be a $1$-cochain of $Y$. We identify $\phi$
with the graph $G=(V,E)$ where $E=\{uv:\phi(uv) =1\}$.  For $v \in V$, let
$$S_v=\Gamma_G(v)=\{u \in V: uv \in E\} \subset V_v$$
and let $s_v=|S_v|=\deg_G(v)$.
\begin{proposition}
\label{onetwo1}
Let $c<1$ and suppose $m=\|[\phi]\| \leq cn^2$. Then
$$
\|d_1\phi\| \geq \left(\frac{(1-c^{1/3})\tilde{\mu}}{2}-\frac{d}{3}\right) \|[\phi]\|.$$
\end{proposition}
\noindent
{\bf Proof:}
We may assume that $\|\phi\|=\|[\phi]\|=m$.
Therefore $\|\phi\| \leq \|\phi+d_0 1_v\|$ and hence $s_v \leq \frac{|V_v|}{2}=n$ for all $v \in V$.
We will need the following
\begin{claim}
\label{ezer2}
$$\sum_{v \in V} s_v^2 \leq (1+3c^{1/3})mn.
$$
\end{claim}
\noindent
{\bf Proof:} Let $$I=\{v \in V: s_v \geq c^{1/3}n \}.$$
Then
$$|I| \leq \frac{\sum_{v \in V} s_v}{c^{1/3}n} =\frac{2m}{c^{1/3}n}.$$
Let
$$
E_0=\{uv \in E: s_u+s_v \geq (1+c^{1/3})n \}.
$$
If $uv \in E_0$ then $u,v \in I$. It follows that
$$|E_0| \leq \frac{|I|^2}{2} \leq
\frac{2m^2}{c^{2/3} n^2} \leq 2 c^{1/3} m.$$
Therefore
\begin{equation*}
\begin{split}
\sum_{v \in V} s_v^2 &= \sum_{uv \in E}(s_u+s_v)= \sum_{uv \in E_0}(s_u+s_v)+ \sum_{uv \in E-E_0}(s_u+s_v) \\
&\leq |E_0|\cdot 2n+(m-|E_0|) \cdot (1+c^{1/3})n  \\
&= (1-c^{1/3})n |E_0|+ (1+c^{1/3}) mn \\
&\leq (1+3c^{1/3})mn.
\end{split}
\end{equation*}
{\enp}
\noindent
{\bf Proof of Proposition \ref{onetwo1}:}
For $v \in V$, define
$\phi_v \in C^{0}(Y;\FF_2)$ by
$$
\phi_v(u)=\left\{
\begin{array}{ll}
\phi(vu) & uv \in Y(1) \\
0 & otherwise
\end{array}
\right.~~
$$
Note that if $uw \in \lk(Y,v)(1)$ then
$$
d_1\phi(vuw)=\phi(uw)-d_{0}\phi_v(uw).
$$
Additionally, if $v \in V$ then
\begin{equation*}
\begin{split}
&|\supp(d_{0}\phi_v) \cap \lk(Y,v)|=|\{uw \in Y_v(1): \phi(vu) \neq \phi (vw)\}| \\
&=|\{uw \in Y_v(1): u \in S_v, w \in V_v-S_v\}|
=e_{Y_v}(S_v,V_v-S_v) .
\end{split}
\end{equation*}
Therefore
\begin{equation*}
\label{xyi}
\begin{split}
3\|d_{1}\phi\|&=\sum_{v \in V}
|\{(v,uw) \in V \times Y(1): uw \in \lk(Y,v)~,~d_{1}\phi(vuw) \neq 0\}| \\
&=\sum_{v \in V} |\supp(\phi-d_{0}\phi_v) \cap \lk(Y,v)| \\
&\geq \sum_{v \in V} |\supp(d_{0}\phi_v) \cap \lk(Y,v)|- \sum_{v \in V} |\supp(\phi) \cap \lk(Y,v)| \\
&= \sum_{v \in V} e_{Y_v}(S_v,V_v-S_v) - \sum_{uw \in \supp(\phi)} |\lk(Y,uw)|.
\end{split}
\end{equation*}
The Alon-Milman and Tanner bound (see Theorem 9.2.1 in \cite{AS}) asserts that
$$e_{Y_v}(S_v,V_v-S_v) \geq \frac{|S_v|(|V_v-S_v|)}{|V_v|} \mu(Y_v) =
\frac{s_v(2n-s_v)}{2n} \mu(Y_v).$$
Combined with the assumption that $|\lk(Y,uw)| \leq d$ for all edges $uw \in T(1)$, it follows that
\begin{equation*}
\begin{split}
3\|d_{1}\phi\| &\geq \tilde{\mu} \sum_{v \in V} \frac{s_v(2n-s_v)}{2n}-md \\
&=\tilde{\mu}(2m-\frac{1}{2n}\sum_{v \in V}s_v^2)-md \\
&\geq \tilde{\mu} (2m-\frac{(1+3c^{1/3})mn}{2n})-md \\
&=\left(\frac{3(1-c^{1/3})\tilde{\mu}}{2}-d\right)m.
\end{split}
\end{equation*}
{\enp}

\section{Large Deviations for Latin Squares}
\label{s:ldls}
\ \ \
Let $0<c \leq 1$ and let $\ce \subset T_n(2)$ be a family of $2$-simplices in $T_n$ such that $|\ce| \geq cn^3$.
Define a random variable $f_{\ce}$ on the space of Latin squares $\cl_n$ by
$$f_{\ce}(L)=|Y(L)(2) \cap \ce|.$$
The expectation of $f_{\ce}$ is
$$E[f_{\ce}]=\frac{|\ce|}{n} \geq cn^2.$$
In the next theorem we establish a one-sided tail estimate for the random variable $f_{\ce}$.
Let us remark that if instead of the space $\cl_n$ we consider the larger probability space $\ccs_n^n$, then a similar estimate follows from Maurey's theorem \cite{Ma79,Milman86}. This however does not seem applicable to our case since the measure of $\cl_n$ inside $\ccs_n^n$ is only $\exp(-\Omega(n^2))$.
\begin{theorem}
\label{largedev}
There exists an $n_0(c)$ such that for all $n \geq n_0(c)$
$$
\pr[f_{\ce}(L)< 10^{-3} c^2 n^2] < e^{-10^{-3} c^2 n^2}.
$$
\end{theorem}

In Subsection \ref{s:lrga} we obtain an upper bound (Proposition \ref{cntb}) on the size of
a certain family of permutations. This is the main ingredient in the proof of Theorem \ref{largedev}
given in Subsection \ref{s:lrgb}.

{\bf
\subsection{Counting Restricted Permutations}
\label{s:lrga}
}
Let $0<\gamma \leq 1/2$ and let $E$  be a subset of $[n]^2$ such that  $|E| \geq \gamma n^2$.
Let $k \leq \frac{\gamma n}{2}$ and let $F=(B_1,\ldots, B_n)$ be an $n$-tuple of subsets of $[n]$ such that $|B_i|=k$ for all $1 \leq i \leq n$.
For a permutation $\pi \in \ccs_n$, let
\begin{equation}
\label{ggee}
g_E(\pi)=|\{(i,\pi(i))\}_{i=1}^n \cap E|.
\end{equation}
Our first goal in this section is to bound the cardinality of the set
$$S(E,F)=\{\pi \in \ccs_n:  g_E(\pi) \leq  \frac{\gamma n}{10}~,~\pi(i) \not\in B_i \text{~for~all~} 1 \leq i \leq n \}.
$$
It is instructive to first consider the case $k=0$ where only the first restriction $g_E(\pi) \leq \frac{\gamma n}{10}$
is relevant. The expectation of $g_E$ over the uniform probability space $\ccs_n$ is $E[g_E]=\frac{|E|}{n} \geq \gamma n$,
hence by Maurey's large deviation theorem \cite{Ma79}
\begin{equation*}
\begin{split}
&\pr\left[\pi \in \ccs_n: g_E(\pi) \leq \frac{\gamma n}{10}\right] \\
&\leq
\pr\left[\pi \in \ccs_n: g_E(\pi) \leq \frac{E[g_E]}{10}\right] < 2 e^{-\frac{\gamma^2 n}{20}}.
\end{split}
\end{equation*}
Therefore if $k=0$ then
\begin{equation}
\label{kzero}
\begin{split}
|S(E,F)| \leq 2 n! \cdot e^{-\frac{\gamma^2 n}{20}} \leq 2e \sqrt{n} \left(\frac{n}{e}\right)^n  e^{-\frac{\gamma^2 n}{20}}.
\end{split}
\end{equation}
Below we prove an extension of (\ref{kzero}) for general $k \leq \frac{\gamma n}{2}$ (under a mild assumption on $E$).
Let us remark that the proof of Maurey's theorem depends on martingale inequalities, while our approach for general $k$ is different.
In fact, the bound we obtain for $k=0$ is slightly better than (\ref{kzero}) when $\gamma$ is small.
\ \\ \\
Write $E=\cup_{i=1}^n (\{i\} \times A_i)$. Then
$$g_E(\pi)=|\{1 \leq i \leq  n: \pi(i) \in A_i\}|.$$
For $0 \leq m \leq n$, let
$$S(E,F,m)=\{\pi \in \ccs_n:  g_E(\pi)=m ~,~\pi(i) \not\in B_i \text{~for~all~} 1 \leq i \leq n \}.
$$
Let
\begin{equation}
\label{elle}
\ell(E)=\max_{1 \leq i \leq n} |A_i|.
\end{equation}
\begin{proposition}
\label{crp}
Let $E,F$ be as above. If $\ell(E)\leq \frac{n}{2}$, then
$$
|S(E,F)|=\sum_{m \leq \frac{\gamma n}{10}} |S(E,F,m)| \leq 4n^2 \left(\frac{n-k}{e}\right)^n e^{-\frac{\gamma n}{20}}.
$$
\end{proposition}
\noindent
{\bf Proof:} For $1 \leq i \leq n$, let
$R_i=A_i-B_i$, $S_i=[n]-A_i-B_i$ and denote $r_i=|R_i|$ and $p_i=\frac{r_i}{n-k}\leq \frac{2}{3}$ (as $r_i \leq \frac{n}{2}$ and $k \leq \frac{n}{4}$).
Then $s_i=|S_i|=n-k-r_i$ and $\frac{s_i}{n-k}=1-p_i$.
For a subset $I \subset [n]$
let $M_I$ be the $n \times n$ matrix given by
$$M_I(i,j)=\left\{
\begin{array}{ll}
1 & i \in I, j \in R_i, \\
1 & i \not\in I, j \in S_i, \\
0 & otherwise.
\end{array}
\right.~~
$$
Note that $\pi \in S(E,F,m)$ iff there exists an $I \in \binom{[n]}{m}$ such that $\pi(i) \in R_i$ for $i \in I$ and
$\pi(i) \in S_i$ for $i \in [n]-I$. The last condition is equivalent to $\prod_{i=1}^n M_I(i,\pi(i))=1$.
As such $I$ must be unique, it follows that
\begin{equation}
\label{sefperm}
|S(E,F,m)|=\sum_{|I|=m} \sum_{\pi \in \ccs_n} \prod_{i=1}^n M_I(i,\pi(i))=\sum_{|I|=m} \per \, M_I.
\end{equation}
Since $s_i \geq n-|A_i|-k \geq \frac{n}{4}$,  the Stirling approximation implies
$$s_i!^{\frac{1}{s_i}} \leq (e^2 s_i)^{\frac{1}{2s_i}} \left(\frac{s_i}{e}\right)  \leq (2n)^{\frac{2}{n}} \left(\frac{s_i}{e}\right).$$
Recall the following result of Br\'{e}gman (see e.g. Theorem 11.5 in \cite{vLW01}).
\begin{theorem}[Br\'{e}gman]
\label{bregman}
Let $M$ be an $n \times n$ zero-one matrix with $t_i$ ones in the $i$-th row. Then
$$
\per \, M \leq \prod_{i=1}^n t_i!^{\frac{1}{t_i}}.
$$
\end{theorem}
Using (\ref{sefperm}) and Br\'{e}gman's bound it follows that
\begin{equation}
\begin{split}
\label{bdd}
|S(E,F,m)|&=\sum_{|I|=m} \per \, M_I \\
&\leq \sum_{|I|=m} \prod_{i \in I} r_i!^{\frac{1}{r_i}} \prod_{i \not\in I} s_i!^{\frac{1}{s_i}} \\
&\leq \sum_{|I|=m} \prod_{i \in I} r_i \prod_{i \not\in I} \left((2n)^{\frac{2}{n}}\left(\frac{s_i}{e}\right) \right) \\
&=(2n)^{\frac{2(n-m)}{n}} e^{m-n} (n-k)^n \sum_{|I|=m} \prod_{i \in I} \left(\frac{r_i}{n-k}\right) \prod_{i \not\in I} \left(\frac{s_i}{n-k}\right) \\
&\leq 4n^2 \cdot  e^m \cdot \left(\frac{n-k}{e}\right)^n \sum_{|I|=m} \prod_{i \in I} p_i \prod_{i \not\in I} (1-p_i).
\end{split}
\end{equation}
Consider a sequence $\{Y_i\}_{i=1}^n$ of independent binomial random variables such that $\pr[Y_i=1]=p_i=1-\pr[Y_i=0]$ and let
$Y=\sum_{i=1}^n Y_i$.
Then by (\ref{bdd})
\begin{equation}
\label{bdda}
|S(E,F,m)| \leq 4n^2 \cdot  e^m \cdot \left(\frac{n-k}{e}\right)^n  \pr[Y=m].
\end{equation}
Next note that
\begin{equation*}
\begin{split}
E[Y]&=\sum_{i=1}^n p_i =\sum_{i=1}^n  \frac{r_i}{n-k} \\
&\geq \frac{1}{n-k}\sum_{i=1}^n(|A_i|-k) \geq \frac{\gamma n^2-kn}{n}
\geq \frac{\gamma n}{2}.
\end{split}
\end{equation*}
The Chernoff bound (see Theorem A.1.13 in \cite{AS}) states that for all $a>0$
\begin{equation}
\label{largedi}
\pr[Y<E[Y]-a]<e^{-\frac{a^2}{2 E[Y]}}.
\end{equation}
Using (\ref{bdda}) and (\ref{largedi}) with $a=\frac{4E[Y]}{5}$ it follows that
\begin{equation}
\label{ld1}
\begin{split}
&(4n^2)^{-1} \cdot  e^{-\frac{\gamma n}{10}} \cdot \left(\frac{n-k}{e}\right)^{-n} \sum_{m< \frac{\gamma n}{10}} |S(E,F,m)| \\
&\leq  \pr \left[ Y<\frac{\gamma n}{10}\right]
\leq  \pr\left[Y<\frac{E[Y]}{5}\right]
\leq  e^{-\frac{4 \gamma n}{25}}
\end{split}
\end{equation}
and therefore
\begin{equation*}
|S(E,F)|=\sum_{m< \frac{\gamma n}{10}} |S(E,F,m)| \leq 4n^2 \cdot \left(\frac{n-k}{e}\right)^n e^{-\frac{\gamma n}{20}}.
\end{equation*}
{\enp}
Now let $E_1,\ldots,E_n$ be subsets of $[n]^2$ such that $\ell(E_i) \leq \frac{n}{2}$ for all $1 \leq i \leq n$.
Suppose $I \subset [n]$ satisfies $|I|=\frac{\gamma n}{2}$ and $|E_i| \geq \gamma n^2$ for all $i \in I$.
We next use Proposition \ref{crp} to bound the number of Latin squares $L=(\pi_1,\ldots,\pi_n)$ such that $g_{E_i}(\pi_i)$ is much smaller than its expected value
for all $i \in I$.
Let
$$\cg(I)= \{L=(\pi_1,\ldots,\pi_n) \in \cl_n: g_{E_i}(\pi_i) < \frac{\gamma n}{10} \text{~for~all~} i \in I\}.$$
The main ingredient in the proof of the large deviation bound for random Latin squares in Section \ref{s:lrgb} is the following
\begin{proposition}
\label{cntb}
$$|\cg(I)| \leq (2n)^{\gamma n} \left(\prod_{k=1}^n k!^{\frac{n}{k}}\right)  e^{-\frac{\gamma^2 n^2}{40}}.$$
\end{proposition}
\noindent
{\bf Proof:} Let $k_0=\frac{\gamma n}{2}$. Without loss of generality we may assume that $I=[k_0]$. A legal $k$-tuple $(\pi_1,\ldots,\pi_k) \in \ccs_n^k$ is {\it extendible to $\cg(I)$} if there exist
$\pi_{k+1},\ldots,\pi_n \in \ccs_n$ such that $(\pi_1,\ldots,\pi_n) \in \cg(I)$.
Fixing a $0 \leq k \leq n-1$ and a legal  $(\pi_1,\ldots,\pi_k) \in \ccs_n^k$,
we next obtain an upper bound on
$$\mu(\pi_1,\ldots,\pi_k)=|\{\pi_{k+1} \in S_k: (\pi_1,\ldots,\pi_{k+1}) ~is~extendible~to~\cg(I)\}|.$$
For $1 \leq i \leq n$, let
 $B_i=\{\pi_j(i):1 \leq j \leq k\}$ and let $F=(B_1,\ldots,B_n)$. If $k<k_0$ then by Proposition \ref{crp}
\begin{equation}
\label{bsmall}
\begin{split}
\mu(\pi_1,\ldots,\pi_k) &\leq \sum_{m< \frac{\gamma n}{10}} |S(E_{k+1},F,m)| \\
&\leq 4n^2 \left(\frac{n-k}{e}\right)^n e^{-\frac{\gamma n}{20}}.
\end{split}
\end{equation}
Suppose on the other hand that $k_0 \leq k <n$. If $(\pi_1,\ldots,\pi_{k+1})$ is extendible to $\cg(I)$
then
$\pi_{k+1}(i) \not\in B_i$ for all $1 \leq i \leq n$. Hence again by Br\'{e}gman's bound
\begin{equation}
\label{bbig}
\mu(\pi_1,\ldots,\pi_k) \leq (n-k)!^{\frac{n}{n-k}}.
\end{equation}
Choosing $\pi_1,\ldots,\pi_n$ sequentially and using (\ref{bsmall}) and (\ref{bbig}) it follows that
\begin{equation*}
\begin{split}
|\cg(I)| &\leq \prod_{k=0}^{k_0-1}\left( 4n^2 \left(\frac{n-k}{e}\right)^n e^{-\frac{\gamma n}{20}}\right)  \prod_{k=k_0}^{n-1}(n-k)!^{\frac{n}{n-k}} \\
&\leq (4n^2 e^{-\frac{\gamma n}{20}})^{\frac{\gamma n}{2}}  \prod_{k=0}^{k_0-1} (n-k)!^{\frac{n}{n-k}} \prod_{k=k_0}^{n-1}(n-k)!^{\frac{n}{n-k}} \\
&=(2n)^{\gamma n} \left(\prod_{k=1}^n k!^{\frac{n}{k}}\right)  e^{-\frac{\gamma^2 n^2}{40}}.
\end{split}
\end{equation*}
{\enp}

{\bf
\subsection{Proof of Theorem \ref{largedev}}
\label{s:lrgb}
}
Recall that $\ce$ is a subset of $2$-simplices of $T_n$ of cardinality $|\ce| \geq cn^3$. We have to bound the probability of
 $$\cl_n(c)=\{L \in \cl_n:  f_{\ce}(L)< 10^{-3} c^2 n^2\}.$$
For $1 \leq k \leq n$, let
$$E_k'=\{(i,j):[a_k,b_i,c_j] \in \ce\} \subset [n]^2$$
and write $E_k'=\cup_{i=1}^n (\{i\} \times J_{ki}')$. For $1\leq k,i \leq n$ choose an arbitrary $J_{ki} \subset J_{ki}'$ such that
$|J_{ki}|=\min\{\frac{n}{2},|J_{ki}'|\}$.
Then $E_k=\cup_{i=1}^n(\{i\} \times J_{ki}) \subset E_k'$ satisfies $\ell(E_k)=\max_{1 \leq i \leq n} |J_{ki}| \leq \frac{n}{2}$
and $|E_k| \geq \frac{|E_k'|}{2}$.
For $L=(\pi_1,\ldots,\pi_n) \in \cl_n$, let
$$g(L)=\sum_{k=1}^n g_{E_k}(\pi_k).$$
Let $\gamma=\frac{c}{4}$. Then
\begin{equation*}
\begin{split}
E[g]&=\frac{1}{n} \sum_{k=1}^n |E_k| \geq \frac{1}{2n} \sum_{k=1}^n |E_k'| \\
&= \frac{|\ce|}{2n} \geq \frac{cn^2}{2}=2\gamma n^2.
\end{split}
\end{equation*}
\begin{claim}
\label{bigset1}
For any $L=(\pi_1,\ldots,\pi_n) \in\cl_n(c)$, there exists a subset $I_L \subset [n]$ of size $|I_L|=\frac{\gamma n}{2}$ such that
for all $i \in I_L$ both $|E_i| \geq \gamma n^2$ and $g_{E_i}(\pi_i) \leq \frac{\gamma n}{10}$.
\end{claim}
\noindent
{\bf Proof:} Let $$I=\{1 \leq i \leq n: |E_i| \geq \gamma n^2\}$$ and
$$J=\{i \in I: g_{E_i}(\pi_i) \leq \frac{\gamma n}{10}\}.$$
Since $\sum_{k=1}^n|E_k| \geq \frac{1}{2}\sum_{k=1}^n |E_k'| \geq 2 \gamma n^3$ it follows that
$|I| \geq \gamma n$.
Therefore
\begin{equation*}
\begin{split}
&(\gamma n-|J|)\frac{\gamma n}{10} \leq (|I|-|J|)\frac{\gamma n}{10}
\leq \sum_{i=1}^n g_{E_i}(\pi_i) \\
&\leq \sum_{i=1}^n g_{E_i'}(\pi_i)=
f_{\ce}(L)<10^{-3} c^2 n^2=16\cdot 10^{-3} \gamma^2 n^2.
\end{split}
\end{equation*}
Hence $|J| \geq (1-\frac{16}{100})\gamma n > \frac{\gamma n}{2}$.  The Claim follows by taking $I_L$ to be any $\frac{\gamma n}{2}$ subset of $J$.
{\enp}
\noindent
Claim \ref{bigset1} implies that $\cl_n(c) \subset \bigcup_{|I|= \frac{\gamma n}{2}} \cg(I)$.
Therefore by Proposition \ref{cntb}
\begin{equation}
\label{clnc}
\begin{split}
|\cl_n(c)| &\leq |\bigcup_{|I|= \frac{\gamma n}{2}} \cg(I)| \\
&\leq \binom{n}{\frac{\gamma n}{2}} (2n)^{\gamma n} \left(\prod_{k=1}^n k!^{\frac{n}{k}}\right)  e^{-\frac{\gamma^2 n^2}{40}} \\
&=\binom{n}{\frac{c n}{8}} (2n)^{\frac{cn}{4}} \left(\prod_{k=1}^n k!^{\frac{n}{k}}\right)  e^{-\frac{c^2 n^2}{640}}.
\end{split}
\end{equation}
By a classical result on the enumeration of Latin squares (See Theorem 17.3  in \cite{vLW01})
\begin{equation}
\label{nls}
1 \leq  \frac{\prod_{k=1}^n  k!^{\frac{n}{k}}}{|\cl_n|}=(1+o(1))^{n^2}.
\end{equation}
It follows by (\ref{clnc}) that for sufficiently large $n \geq n_0(c)$
$$
\pr[\cl_n(c)] =\frac{|\cl_n(c)|}{|\cl_n|} \leq
e^{-10^{-3} c^2 n^2}.
$$
{\enp}

\section{Expanders from Latin Squares}
\label{prf}
\ \ \
Here we prove Theorem \ref{main}. We first consider the expansion of small cochains.
For a $d$-tuple of permutations $\tilde{\pi}=(\pi_1,\ldots,\pi_d) \in \ccs_n^d$, let
$G=G(\tilde{\pi})$ be the $n$ by $n$ bipartite graph whose edge set
is $$\{(i,\pi_j(i)):1 \leq i \leq n, 1 \leq j \leq d\} \subset [n]^2.$$
Let $\cg(n,d)$ be the uniform probability space $\{G(\tilde{\pi}):\tilde{\pi} \in \ccs_n^d \}$.
Friedman's theorem \cite{F91,F08} on the spectral gap of $d$-regular graphs implies that if $d \geq 100$ is fixed then
$\mu(G)>d-3\sqrt{d}$ with with probability $1-O(n^{-2})$.

Consider a random $Y(\ul^d) \in \cy(n,d)$.
The link $\lk(Y(\ul^d),v)$ of a fixed vertex $v \in V$ is a random graph in $\cg(n,d)$ hence
$$\pr[\mu(\lk(Y(\ul^d),v))]>d-3\sqrt{d}]=1-O(n^{-2}).$$ Since $|V|=3n$ it follows that
$$\pr[\min_{v \in V} \mu(\lk(Y(\ul^d),v))>d-3\sqrt{d}]=1-O(n^{-1}).$$
Proposition \ref{onetwo1} thus implies the following
\begin{corollary}
\label{corr1}
For any fixed $d \geq 100$ and $c<1$, the following holds with probability $1-O(n^{-1})$:
$$\frac{\|d_1\phi\|_{Y(\ul^d)}}{\|[\phi]\|} \geq \frac{(d-3\sqrt{d})(1-c^{1/3})}{2}-\frac{d}{3}$$
for all $1$-cochains $\phi \in C^1(T_n;\FF_2)$ such that $\|[\phi]\| \leq cn^2$.
{\enp}
\end{corollary}

We next consider the expansion of large cochains.
Fix $\phi \in C^1(T_n;\FF_2)$ such that $\|[\phi]\| \geq cn^2$ with $c<1$.
A special case of a result of Dotterrer and Kahle
(Proposition 5.7 in \cite{DK12}) asserts that $h_1(T_n) \geq \frac{n}{5}$.
It follows that
 $$\ce=\{\sigma \in T_n(2): d_1\phi(\sigma) \neq 0\}$$
satisfies
$$|\ce|=\|d_1\phi\|_{T_n} \geq \frac{n}{5} \|[\phi]\| \geq \frac{cn^3}{5}.$$
If $L \in \cl_n$ then $\|d_1\phi\|_{Y(L)}=|Y(L)(2) \cap \ce|=f_{\ce}(L)$.
Theorem \ref{largedev} then implies that if $n \geq n_0(c/5)$ then
\begin{equation}
\label{pln}
\begin{split}
& \pr[L \in \cl_n: \|d_1\phi\|_{Y(L)} < 10^{-3}(c/5)^2n^2] \\
&= \pr[L \in \cl_n: f_{\ce}(L) < 10^{-3}(c/5)^2n^2]  \\
&< e^{-10^{-3}(c/5)^2n^2}=e^{-4\cdot 10^{-5}c^2n^2}.
\end{split}
\end{equation}
Let $\ul^d=(L_1,\ldots,L_d) \in \cl_n^d$. Then
$$\|d_1\phi\|_{Y(\ul^d)} \geq \max_{1 \leq i \leq d} \|d_1\phi\|_{Y(L_i)}.$$
Therefore by (\ref{pln})
$$\pr[\ul^d \in \cl_n^d: \|d_1\phi\|_{Y(\ul^d)} < 4 \cdot10^{-5}c^2n^2] <e^{-4\cdot 10^{-5}dc^2n^2}.$$
As the number of $1$-cochains is $2^{3n^2}$ and the support of a $1$-cochain is at most $3n^2$ we obtain
\begin{corollary}
\label{corr2}
If $n \geq n_0(c/10)$ then the following holds with probability at least  $1-2^{3n^2}e^{-4\cdot 10^{-5}dc^2n^2}$:
$$\frac{\|d_1\phi\|_{Y(\ul^d)}}{\|[\phi]\|} \geq 10^{-5}c^2$$
for all $1$-cochains $\phi \in C^1(T_n;\FF_2)$ such that $\|\phi\| \geq cn^2$.
{\enp}
\end{corollary}
\ \\ \\
{\bf Proof of Theorem \ref{main}:} Let $c=10^{-3}$ then
$$\frac{(d-3\sqrt{d})(1-c^{1/3})}{2}-\frac{d}{3}>1$$
for $d>200$.
Let $d=10^{11}$ then
$2^3<e^{4\cdot 10^{-5}dc^2}$
so Corollaries \ref{corr1} and \ref{corr2}
imply that Theorem \ref{main} holds with $d=10^{11}$ and $\epsilon=10^{-5}c^2=10^{-11}$.
{\enp}

\section{Concluding Remarks}
\label{crem}
\ \ \
We have shown that there exist fixed $d <\infty$ and $\epsilon>0$ such that $Y \in \cy(n,d)$ a.a.s. satisfies $h_1(Y)> \epsilon$. In particular, there exist infinite families of $(2,d,\epsilon)$-expanders.
We conclude with the following comments and questions.
\begin{enumerate}
\item
Let $d_0$ be the smallest $d$ for which Theorem \ref{main} remains true.
It can be shown that if $d=3$ then
$$\lim_{n \rightarrow \infty} \pr[H_1(Y(L_1,L_2,L_3);\FF_2) \neq 0] \geq 1-\frac{17 e^{-3}}{2}\doteq 0.57 .$$
In particular $d_0 \geq 4$.  It seems plausible that $d_0$ is in fact $4$. A considerably weaker question would be to determine for a fixed field $\FF$ the minimal $d$ such that $H_1(Y;\FF)=0$ a.a.s. for $Y \in \cy(n,d)$.
For $\FF=\Rea$ one can use a spectral approach similar to the one applied in \cite{GW12} and \cite{HKP12} for other models.
A classical result of Garland \cite{Gar73} asserts (roughly) that if in a $2$-dimensional complex $Y$ all vertex links have sufficiently large spectral gaps then $H_1(Y;\Rea)=0$. Combining Friedman's eigenvalue bounds for random graphs
\cite{F91,F08} and Garland's theorem it follows that if $d \geq 100$ then $H_1(Y;\Rea)=0$ a.a.s. for
$Y \in \cy(n,d)$.

\item
The complexes $Y \in \cy(n,d)$ satisfy $D_1(Y) \leq d$ but $D_0(Y) \geq n$.
It would be very interesting to prove the existence (or even better to give explicit constructions) of infinite families of $\epsilon$-expanding complexes $X$ such that both $D_0(X)$ and $D_1(X)$ are bounded. For some results in this direction see the work of Kaufman, Kazhdan and Lubotzky \cite{KKL13} on expansion in Ramanujan complexes.

\item
The model $\cy(n,d)$ generalizes in a straightforward manner to higher dimensions and it seems likely that Theorem \ref{main} remains true there.
The main obstacle to extending the present proof to this case is the absence (at present) of higher dimensional analogues of the asymptotic enumeration formula
(\ref{nls}).

\item
A simplicial complex $X$ is said to have the $(k,\delta)$ topological overlap property if for any continuous map $f:X \rightarrow \Rea^k$ there exists a point $p \in \Rea^k$
such that
$$|\{\sigma \in X(k): p \in f(\sigma)\}| \geq \delta |X(k)|.$$
A remarkable theorem of Gromov \cite{G10} asserts that for any $k$ and $\epsilon>0$ there exists an $\delta=\delta(k,\epsilon)$ such that if
$h_i(X) \geq \epsilon \cdot \frac{|X(i+1)|}{|X(i)|}$
for every $0 \leq i \leq k-1$, then $X$ has the $(k,\delta)$ topological overlap property.
Theorem \ref{main} therefore implies that there exist $d$ and $\delta>0$ such that $Y \in Y(n,d)$ a.a.s. has the $(2,\delta)$ topological overlap property.

\end{enumerate}


\begin{thebibliography}{99}

\bibitem{AS}
N. Alon and J. Spencer, {\it The Probabilistic Method}, 3rd
Edition, Wiley-Intescience, 2008.

\bibitem{DK12}
D. Dotterrer and M. Kahle, Coboundary expanders, {\it J. Topol. Anal.} {\bf 4}(2012) 499-–514.


\bibitem{FGLNP}
J. Fox, M. Gromov, V. Lafforgue, A. Naor and J. Pach, Overlap properties of geometric expanders, {\it J. Reine Angew. Math.} {\bf 671}(2012) 49–-83.

\bibitem{F91}
J. Friedman, On the second eigenvalue and random walks in random d-regular graphs, {\it Combinatorica} {\bf 11} (1991) 331-–362.

\bibitem{F08}
J. Friedman, A proof of Alon's second eigenvalue conjecture and related problems. Mem. Amer. Math. Soc.  {\bf 195}(2008).

\bibitem{Gar73}
H. Garland, $p$-adic curvature and the cohomology of discrete
subgroups of $p$-adic groups, {\it Annals of Math.} {\bf 97}
(1973) 375-423.

\bibitem{G10}
M. Gromov, Singularities, expanders and topology of maps. Part 2: From combinatorics to topology via algebraic isoperimetry, {\it Geom. Funct. Anal.} {\bf 20}(2010) 416-–526.

\bibitem{GW12}
A. Gundert and U. Wagner,
On Laplacians of Random Complexes,
Proc. 28th Annual ACM Symposium on Computational Geometry (2012) 151--160.

\bibitem{HKP12}
C. Hoffman, M. Kahle and E. Paquette,
A sharp threshold for Kazhdan's property (T),
arXiv:1201.0425

\bibitem{HLW06}
S. Hoory, N. Linial and A. Wigderson, Expander graphs and their applications, {\it Bull. Amer. Math. Soc. (N.S.)} {\bf 43}(2006) 439–-561.


\bibitem{KKL13}
T. Kaufman, D. Kazhdan and A. Lubotzky, Isoperimetric inequalities for Ramanujan complexes, in preparation.

\bibitem{LM06}
N. Linial and R. Meshulam, Homological connectivity of random
2-complexes, {\it Combinatorica} {\bf 26}(2006) 475--487.

\bibitem{vLW01}
J. H. van Lint and R. M. Wilson, A Course in Combinatorics. Second
edition. Cambridge University Press, Cambridge, 2001.

\bibitem{Lu1}
A. Lubotzky, Discrete Groups, Expanding Graphs and Invariant Measures. With an appendix by Jonathan D. Rogawski. Progress in Mathematics, 125. Birkh\"{a}user Verlag, Basel, 1994.

\bibitem{Lu2}
A. Lubotzky, Expander graphs in pure and applied mathematics, {\it Bull. Amer. Math. Soc. (N.S.)} {\bf 49}(2012) 113–-162.

\bibitem{Lu3}
A. Lubotzky, Ramanujan complexes and high dimensional expanders, Takagi Lectures, Tokyo, 2012,
arXiv:1301.1028.

\bibitem{Ma79}
B. Maurey, Construction de suites sym\'{e}triques, C. R. Acad. Sci. Paris S\'{e}r. A-B {\bf 288}(1979), no. 14, A679-–A681.

\bibitem{MW09}
R. Meshulam and N. Wallach, Homological connectivity of
random $k$-dimensional complexes,
{\it Random Struct. Algorithms}  {\bf 34}(2009) 408--417.

\bibitem{Milman86}
V. D. Milman, G. Schechtman, Asymptotic Theory of Finite-Dimensional Normed Spaces. With an appendix by M. Gromov. Lecture Notes in Mathematics, 1200. Springer-Verlag, Berlin, 1986.

\end{thebibliography}
\end{document}